\documentclass{article}
\usepackage[utf8]{inputenc}
\usepackage[T1]{fontenc}

\usepackage{amssymb}
\usepackage{stmaryrd}
\usepackage{amsmath,amsfonts}
\usepackage[tikz]{bclogo}
\usepackage[top=4.34cm]{geometry}

\usepackage{lmodern}
\usepackage{textcomp}
\usepackage{mathtools, bm}
\usepackage{amssymb, bm}
\usepackage[unq]{unique}
\usepackage{enumerate}
\usepackage{bbm}

\usepackage{algorithm}
\usepackage{algpseudocode}

\usepackage[breaklinks]{hyperref}
\usepackage{cleveref}
\usepackage[numbers]{natbib}    

\usepackage{amsthm}
\usepackage{thmtools}
\usepackage{comment}
\usepackage{todonotes}

\newtheorem{theorem}{Theorem}[section]
\newtheorem{lemma}[theorem]{Lemma}

\newtheorem{proposition}[theorem]{Proposition}

\theoremstyle{definition}

\title{Ternary Egyptian fractions with prime denominator}
\author{Adva Mond and Julien Portier \thanks{Department of Pure Mathematics and Mathematical Statistics (DPMMS), University of Cambridge, Wilberforce Road, Cambridge, CB3 0WA, United Kingdom. Email: \href{mailto:jp899@cam.ac.uk}{\nolinkurl{{am2759,jp899}@cam.ac.uk}}}}
\date{}

\begin{document}

\maketitle

\begin{abstract}
For a prime number $p$, let $A_3(p)= | \{ m \in \mathbb{N}: \exists m_1,m_2,m_3 \in \mathbb{N}, \frac{m}{p}=\frac{1}{m_1}+\frac{1}{m_2}+\frac{1}{m_3} \} |$.
In 2019 Luca and Pappalardi proved that $x (\log x)^3 \ll \sum_{p \le x} A_{3}(p) \ll x (\log x)^5$.
We improve the upper bound, showing $\sum_{p \le x} A_{3}(p) \ll x (\log x)^3 (\log \log x)^2$.
    
\end{abstract}

\section{Introduction}
An \emph{Egyptian fraction} is a representation of a rational number as a sum of reciprocals of distinct integers.
A ternary Egyptian fraction is such a sum that consists of exactly three summands.
More precisely, it is a representation of a rational number $\frac{m}{n}$ as the sum $\frac{m}{n} = \frac{1}{m_1} + \frac{1}{m_2} + \frac{1}{m_3}$, for some distinct integers $m_1, m_2, m_3$.

Questions regarding Egyptian fractions are amongst the most ancient problems in mathematics.
Throughout history many mathematicians have studied this topic, gaining popularity in recent times thanks to Erd\H{o}s who presented and solved various problems concerning Egyptian fractions (for more details, see, e.g.~\cite{Guy}).
Probably one of the most famous amongst them is a conjecture by Erd\H{o}s and Straus, stating that for any $n \ge 2$, the rational number $\frac{4}{n}$ has a representation as a ternary Egyptian fraction, that is, that the Diophantic equation
\[\frac{4}{n} = \frac{1}{x} + \frac{1}{y} + \frac{1}{z} \]
has at least one solution.
This conjecture is still open.

In this paper we consider ternary Egyptian fractions for which the denominator is a prime number, and we are interested in bounding the number of those, for all primes in a certain range.
As usual, for two functions $f, g : \mathbb N \rightarrow \mathbb R$, by $f(x) \ll g(x)$ we mean that there exists a constant $c > 0$ and a natural number $N \in \mathbb N$, such that for any $n \ge N$ we have $f(n) \le c \cdot g(n)$.
Throughout the paper, $p$ always designates a prime number. \\

Let $A_3(p)= \left| \left\{ m \in \mathbb{N}: \exists m_1,m_2,m_3 \in \mathbb{N}, \frac{m}{p}=\frac{1}{m_1}+\frac{1}{m_2}+\frac{1}{m_3} \right\} \right|$.
Luca and Pappalardi~\cite{LucaPappalardi} proved the following.
\begin{theorem}
\label{thm:LucaPappalardi}
    \[x (\log x)^3 \ll \sum_{p \leq x} A_{3}(p) \ll x (\log x)^5. \]
\end{theorem}

Our main result in this paper closes the gap between the upper and the lower bounds up to a factor of polyloglog.

\begin{theorem}
\label{thm:main}
\begin{align}
\label{eq:result}
    x (\log x)^3 \ll \sum_{p \leq x} A_{3}(p) \ll x (\log x)^3 (\log \log x)^2.
\end{align}
\end{theorem}

Throughout the paper, $\log$ always stands for the logarithm function in base $2$.

\section{Proof idea}
The proof of our main theorem follows the lines of the proof of \Cref{thm:LucaPappalardi} by Luca and Pappalardi~\cite{LucaPappalardi}.
Our contribution is the improved upper bound in \Cref{lem:main}, which is our main lemma.
The proof of \Cref{lem:main} is based on two ingredients.
The first one is an application of the Brun-Titchmarsh inequality (\Cref{brunTitchmarsh}).
The second ingredient is \Cref{prop:ElsholtzTaoImproved}, which is a strengthened version of \Cref{prop:ElsholtzTao} for the certain range of parameters which fits our needs.\\

The next lemma describes a well-known classification of solutions for ternary Egyptian fractions with a prime denominator.
It appears in Mordell's book~\cite{Mordell}, for example, as well as in other texts.
A proof can be found, e.g., in~\cite{LucaPappalardi}.

\begin{lemma}
\label{lem:types}
	If $\frac{m}{p}=\frac{1}{m_1}+\frac{1}{m_2}+\frac{1}{m_3}$, where $m_1, m_2, m_3$ are positive integers and $\gcd (m,p)=1$, then either $m \in \{ 1,2,3 \}$ or there exist positive integers $a, b, c, u$ such that $\gcd (a,b)=1$, $c | a+b$ and one of the following holds:
	
	\begin{itemize}
		\item either (Type I)
		\[m = \frac{p + (a+b)/c}{abu}, \]
		
		\item or (Type II)
		
		\[m = \frac{1 + p(a+b)/c}{abu}. \]
	\end{itemize}
\end{lemma}

Given \Cref{lem:types}, we denote by $A_{3, I}(p)$ and by $A_{3, II}(p)$ the number of those $m\in \mathbb N$ for which $\frac{m}{p}$ is of type I and of type II, respectively.
Given the lower bound of \Cref{thm:LucaPappalardi}, we can already rule out the case where $m \in \{1,2,3 \}$ or $\gcd(m,p) > 1$, as it contributes $O(x)$ to the sum below.
Hence,
\begin{align*}
    \sum_{p\le x} A_3(p) \ll \sum_{p\le x} A_{3,I}(p) + \sum_{p\le x} A_{3,II}(p).
\end{align*}

In \cite{LucaPappalardi}, they deduce \Cref{thm:LucaPappalardi} from the following lemma.

\begin{lemma}
\label{lem:LucaPappalardi}
	We have $x (\log x)^3 \ll \sum_{p \leq x} A_{3,I}(p) \ll x (\log x)^3$ and $\sum_{p \leq x} A_{3,II}(p) \ll x (\log x)^5$.
\end{lemma}

We improve the upper bound on the sum of solutions of type II in \Cref{lem:LucaPappalardi}.

\begin{lemma}
\label{lem:main}
	We have
	\begin{align*}
	    \sum_{p \leq x} A_{3,II}(p) \ll x (\log x)^3 (\log \log x)^2.
	\end{align*}
\end{lemma}

\Cref{thm:main} then follows immediately from \Cref{lem:main}.
The rest of the paper is dedicated to proving \Cref{lem:main}.

\section[Proof of Lemma 2.3]{Proof of \Cref{lem:main}}
We use two classical number theory inequalities.
The first one is the Brun-Titchmarsh inequality (Theorem 6.6 in \cite{AnalyticNT}).
Let $\pi(x;q,a)$ denote the number of primes $p$ congruent to $a$ modulo $q$ satisfying $p \le x$.
Recall that $\phi$ is the Euler totient function.

\begin{theorem}
\label{brunTitchmarsh}
    For all $q < x$ we have
    \begin{align*}
        \pi(x;q,a) \leq \frac{2x}{\phi(q) \log(x/q)}.
    \end{align*}
\end{theorem}

The second inequality we use is the known bound on the sum of characters by Burgess~\cite{burgess}.

\begin{theorem}
\label{Burgess}
    Let $\chi$ be a Dirichlet character modulo $q$.
    Let $r \ge 1$, $H \ge 1$ be fixed integers, and fix $\varepsilon > 0$.
    Then if either $q$ is square-free or $r=2$ we have
    \begin{align*}
        \sum_{N \le n \le N+H} \chi(n) \ll_{r, \varepsilon} H^{1-\frac{1}{r+1}} q^{\frac{1}{4r} + \varepsilon}.
    \end{align*}
\end{theorem}

Recall that $\tau(n) \coloneqq \sum_{d|n} 1$ is the number of distinct divisors $d$ of $n$. Elsholtz and Tao proved the following (Proposition 1.4 from \cite{ElsholtzTao}).

\begin{proposition}
\label{prop:ElsholtzTao}
    For any $A,B > 1$, and any positive integer $k \leq (AB)^{O(1)}$, we have 
    \begin{align*}
        \sum_{a \leq A} \sum_{b \leq B} \tau(kab^2+1) \ll AB \log(A+B) \log(1+k).
    \end{align*}
\end{proposition}

For our proof we need a refined version of \Cref{prop:ElsholtzTao}, which holds for a more restricted range of $k$.
\begin{proposition}
\label{prop:ElsholtzTaoImproved}
For any $A,B > 1$ and $p < \frac{5}{3}$, and any positive integer $k \leq A^p$, we have
\begin{align}
\label{eq:ElsholtzTaoImproved}
    \sum_{a \leq A} \sum_{b \leq B} \tau(kab^2+1) \ll_p AB \log(A+B).
\end{align}
\end{proposition}

The tighter upper bound of \Cref{prop:ElsholtzTaoImproved} is one of the main ingredients in our improved upper bound in \Cref{lem:main}.
Note that \Cref{prop:ElsholtzTaoImproved} can probably be proved for a larger range than $k \leq A^p$ for $p < \frac{5}{3}$, but since in our proof we use \Cref{prop:ElsholtzTaoImproved} only for $p=1$, we have not made any effort in this direction.

\begin{proof}
The proof follows the same lines as of the proof of \Cref{prop:ElsholtzTao} by Elsholtz and Tao~\cite{ElsholtzTao}.
For the case $A \ge B$ it was already shown in~\cite{ElsholtzTao} that (\ref{eq:ElsholtzTaoImproved}) holds.
For the case where $A \le B$, using the same argument as in their proof, it is sufficient to show that
\begin{align*}
    \left|\sum_{\substack{q \le B, \\ (q,2k)=1}} \sum_{\substack{a \le A, \\ (a,2q)=1}} \left( \frac{-ka}{q} \right ) \frac{\log \left(\frac{B}{q} \right)}{q} \right| \ll_p A \log B,
\end{align*}
where by $\left(\frac{a}{q} \right)$ we mean the Jacobi symbol.
Moreover, the contribution of $q > kA$ has been shown by Elsholtz and Tao to be at most $A\log B$.

It is left to consider the contribution of $q \le kA$, for which we obtain a stronger upper bound, using \Cref{Burgess} for $r=2$, and $k \le A^p$.
Thus, we have
\begin{align*}
    \left| \sum_{\substack{a \le A, \\ (a,2q)=1}} \left(\frac{-ka}{q} \right) \right| \ll_{\varepsilon} A^{\frac{2}{3}}q^{\frac{1}{8}+\varepsilon}
\end{align*}

Hence,
\begin{align*}
    \left|\sum_{\substack{q \le kA, \\ (q,2k)=1}} \sum_{\substack{a \le A, \\ (a,2q)=1}} \left( \frac{-ka}{q} \right) \frac{\log \left(\frac{B}{q} \right)}{q} \right| &\ll_{\varepsilon} \sum_{q \leq kA} A^{\frac{2}{3}}q^{\frac{1}{8}-1+\varepsilon} \log B\\
    &\ll_{\varepsilon} A^{\frac{2}{3}}(kA)^{\frac{1}{8}+\varepsilon} \log B \\
    &\ll_{\varepsilon} A^{\frac{2}{3}+\frac{p+1}{8}+\varepsilon(p+1)} \log B,
\end{align*}
Taking $\varepsilon > 0$ small enough proves the statement.
\end{proof}

We are now ready to prove our main lemma.

\begin{proof}[Proof of \Cref{lem:main}.]
In fact, we prove something slightly stronger.
We bound the number of tuples $(m,p,a,b,c,u)$ of Type II satisfying \Cref{lem:types}, which we denote by $\mathcal T(x)$.
This gives an upper bound on the number of pairs $(m,p)$ of type II,
\begin{align}
\label{eq:Tlwr}
    \sum_{p \le x} A_{3,II}(p) \le \mathcal T(x).
\end{align}

For each pair $(m,p)$ of type II we can write $p = \frac{baum - 1}{(a+b)/c}$.
By setting $t = (a+b)/c$ and substituting $b = ct-a$, we get
\begin{align*}
    p = \frac{(ct-a)aum - 1}{t} = caum - \frac{a^2um+1}{t}.
\end{align*}

Furthermore, note that $aum \le 4x$.
Indeed, assuming without loss of generality that $a \le b$, we get that
\begin{align*}
    m = \frac{1+pt}{a(ct-a)u} \le \frac{2pt}{a(ct/2)u} = \frac{4p}{acu},
\end{align*}
giving $aum \le \frac{4p}{c} \le 4x$.
For the sake of simplicity, as $aum \ll x$, we might as well assume $aum \le x$.
Moreover, by symmetry, we can assume $u \le m$.
We have $\tau(a^2um+1)$ possibilities for $t$, and once $a$, $u$, $m$ and $t$ have been fixed, there are only $\pi(x;aum,-(a^2um + 1)/t)$ possibilities for $p$.
Hence, the number of tuples $(m,p,a,b,c,u)$ is at most
\begin{align}
\label{eq:Tupr}
    \mathcal T(x) \le \sum_{aum \le x}\sum_{t | a^2um +1} \pi \left(x; aum, -\frac{a^2um + 1}{t} \right).
\end{align}
Considering both (\ref{eq:Tlwr}) and (\ref{eq:Tupr}), we now focus on bounding from above the right-hand side of (\ref{eq:Tupr}).

By the Brun-Titchmarsh inequality (\Cref{brunTitchmarsh}), we have $\pi(x; aum, d) \ll \frac{x}{\phi(aum)\log(x/aum)}$ for all $d$.
Moreover, considering also the trivial bound $\pi(x; aum, d) \ll \frac{x}{aum}$, we actually get $\pi(x; aum, d) \ll \frac{x}{\phi(aum) \log(2+x/aum)}$, which is useful for those values of $aum$ which are very close to $x$.
Using this last inequality and the classical inequality $\phi(n) \gg \frac{n}{\log \log n}$, we have

\begin{align*}
    \sum_{aum \le x}\sum_{t | a^2um +1} \pi \left(x; aum, -\frac{a^2um + 1}{t} \right) &\ll
    \sum_{aum \le x} \sum_{t|a^2um+1} \frac{x}{\phi(aum) \log(2+x/aum)} \\
    &\ll \sum_{aum \le x} \tau(a^2um+1)\frac{x}{\phi(aum) \log(2+x/aum)} \\
    &\ll \sum_{aum \le x} \tau(a^2um+1)\frac{x}{aum} \frac{\log \log aum}{\log(2+x/aum)}.
\end{align*}
	
It suffices to show that the following holds for any $N \le x$,
\begin{align}
\label{eq:dyadicsum}
    \sum_{N/2 \leq aum \leq N} \frac{\tau(a^2um+1)}{aum} \ll (\log x)^3.
\end{align}
Indeed, summing (\ref{eq:dyadicsum}) over all $N = 2^i$ for $i \leq \log x$ gives
\begin{align*}
    \sum_{p \leq x} A_{3,II}(p) \ll x (\log x)^3 \sum_{i=1}^{\log x} \frac{\log i}{1+\log x -i} \ll x (\log x)^3 (\log \log x)^2,
\end{align*}
proving the lemma. \\

Hence, it is left to prove (\ref{eq:dyadicsum}).
We have
\begin{align*}
    \sum_{N/2 \le aum \le N} \frac{\tau(a^2um+1)}{aum} \ll \sum_{A,U,M} \sum_{U \le u \le 2U} \sum_{A \le a \le 2A} \sum_{M \le m \le 2M} \frac{\tau(a^2um+1)}{aum},
\end{align*}
where the first sum on the right-hand-side is going over all dyadic triplets $(A,U,M) = (2^i, 2^j, 2^h)$ for which the set $\left\{aum ~:~ A\le a\le 2A, \, U\le u \le 2U, \, M\le m \le 2M \right\}$ has a non-empty intersection with the interval $[N/2, N]$.

By \Cref{prop:ElsholtzTaoImproved}, since $U \le M$, we have
\begin{align*}
    \sum_{U \leq u \leq 2U} \sum_{A \leq a \leq 2A} \sum_{M \leq m \leq 2M} \tau(a^2um+1) \ll AUM \log x.
\end{align*}
Since in this range of summation we have $aum \ge AUM$, we get
\begin{align}
\label{eq:aumsum}
    \sum_{U \leq u \leq 2U} \sum_{A \leq a \leq 2A} \sum_{M \leq m \leq 2M} \frac{\tau(a^2um+1)}{aum} \ll \log x.
\end{align}

For every $N \le x$ there are $O((\log x)^2)$ dyadic triplets $(A,U,M)$ for which the set $\{ aum: A \leq a \leq 2A, \, U \leq u \leq 2U, \, M \leq m \leq 2M \}$ has a non-empty intersection with $[N/2,N]$.
Considering (\ref{eq:aumsum}) we then get
\begin{align*}
    \sum_{N/2 \leq aum \leq N} \frac{\tau(a^2um+1)}{aum} \ll (\log x)^3,
\end{align*}
proving (\ref{eq:dyadicsum}), as desired.
\end{proof}

\section{Concluding remarks}

We believe that the correct order is the lower bound $x(\log x)^3$. As mentioned at the beginning of the proof of \Cref{lem:main}, we actually count tuples $(m,p,a,b,c,u)$ rather than pairs $(m,p)$.
A more direct count of the number of pairs $(m,p)$ could possibly yield the desired order of $x(\log x)^3$.

\section{Data availability statement}

Data sharing not applicable to this article as no datasets were generated or analysed during the current study.

\section*{Acknowledgement}

The authors would like to thank their PhD supervisor Professor Béla Bollobás for his valuable comments.\\

In a previous version of this paper we proved an upper bound of $x (\log x)^3 (\log \log x)^3$. We would like to thank Matteo Bordignon, Christian Elsholtz, Bryce Kerr and Timothy Trudgian for pointing out to us that using the Burgess bound instead of the P\'{o}lya-Vinogradov inequality enables us to prove \Cref{prop:ElsholtzTaoImproved} in its current more general version, and consequently removes one $\log \log x$ factor in \Cref{lem:main}.
The authors would also like to thank the anonymous referee for further comments.

\bibliographystyle{abbrvnat}  
\renewcommand{\bibname}{Bibliography}
\bibliography{bibliography}

\begin{thebibliography}{6}
\providecommand{\natexlab}[1]{#1}
\providecommand{\url}[1]{\texttt{#1}}
\expandafter\ifx\csname urlstyle\endcsname\relax
  \providecommand{\doi}[1]{doi: #1}\else
  \providecommand{\doi}{doi: \begingroup \urlstyle{rm}\Url}\fi

\bibitem[Burgess(1963)]{burgess}
D.~A. Burgess.
\newblock On character sums and l-series. ii.
\newblock \emph{Proceedings of the London Mathematical Society}, 3\penalty0
  (1):\penalty0 524--536, 1963.

\bibitem[Elsholtz and Tao(2013)]{ElsholtzTao}
C.~Elsholtz and T.~Tao.
\newblock Counting the number of solutions to the {E}rdős-{S}traus equation on
  unit fractions.
\newblock \emph{Journal of the Australian Mathematical Society}, 94:\penalty0
  50--105, 2013.

\bibitem[Guy(1994)]{Guy}
R.~Guy.
\newblock \emph{Unsolved Problems in {N}umber {T}heory}.
\newblock New York: Springer-Verlag, 2nd edition, 1994.

\bibitem[Iwaniec and Kowalski(2004)]{AnalyticNT}
H.~Iwaniec and E.~Kowalski.
\newblock Analytic number theory.
\newblock \emph{American Mathematical Society Colloquium Publications}, 53,
  2004.

\bibitem[Luca and Pappalardi(2019)]{LucaPappalardi}
F.~Luca and F.~Pappalardi.
\newblock On ternary egyptian fractions with prime denominator.
\newblock \emph{Research in Number Theory}, 5\penalty0 (4):\penalty0 1--14,
  2019.

\bibitem[Mordell(1969)]{Mordell}
L.~J. Mordell.
\newblock \emph{Diophantine equations}.
\newblock Academic Press, 1969.

\end{thebibliography}

\end{document}